\documentclass[11pt]{amsart}
\usepackage{amstext,amssymb,amsmath,amsbsy,dsfont,tikz}

 \usepackage[left=2.7cm,right=2.7cm,top=3.5cm,bottom=3.5cm]{geometry}
 
\usepackage{amsmath,amssymb,latexsym,dsfont}
\usepackage[small]{caption}
\usepackage{graphicx,color,mathrsfs,tikz}%wasysym,overpic,tikz,
\usepackage{subfigure,color}
\usepackage{cite}
\usepackage[colorlinks=true,urlcolor=blue,
citecolor=red,linkcolor=blue,linktocpage,pdfpagelabels,
bookmarksnumbered,bookmarksopen]{hyperref}
\usepackage[italian,english]{babel}
\usepackage{units}
\usepackage{enumitem}
  
\usepackage{hyperref}
\usepackage{cleveref}

\usepackage{tikz}
\usetikzlibrary{intersections}

\usepackage{amscd}
\usepackage{amsfonts}
\usepackage{indentfirst}
\usepackage{verbatim}
\usepackage{amsmath}
\usepackage{amsthm}
\usepackage{enumerate}
\usepackage{graphicx}
\usepackage{color}
\usepackage[OT1]{fontenc}
\usepackage[latin1]{inputenc}
\usepackage[english]{babel}
\usepackage{amssymb,esint}
\newtheorem{theorem}{Theorem}[section]
\newtheorem{lemma}{Lemma}[section]
\newtheorem{proposition}{Proposition}[section]

\newtheorem{remark}{Remark}[section]

\usepackage{tikz}
\usetikzlibrary{intersections}

\numberwithin{equation}{section}

\newcommand{\tr}{^\mathsf{T}}

\newcommand{\hgamma}{\hat \gamma}

\newcommand{\M}{{\mathcal M}}

\newcommand{\supp}{\operatorname{supp}}
\newcommand{\dist}{\operatorname{dist}}

\newcommand{\dive}{\operatorname{div}}

\newcommand{\eps}{\varepsilon}

\newcommand{\mR}{\mathbb{R}}

\newcommand{\hvarphi}{\hat \varphi}
\newcommand{\hLambda}{\hat \Lambda}
\newcommand{\hM}{\hat M}

\newcommand{\cR}{{\mathcal R}}

\newcommand{\tx}{\widetilde x}

\newcommand{\rr}{\hat r}

\title[Cost of fast controls  for the heat equation]{A dependence of the cost of fast controls  for the heat equation on the support of initial datum}

\author{Hoai-Minh Nguyen}
\address[Hoai-Minh Nguyen]{Ecole Polytechnique F\'ed\'erale de Lausanne, EPFL,  CAMA, Station 8,  CH-1015 Lausanne, Switzerland.}
\email{hoai-minh.nguyen@epfl.ch}

\begin{document}
\maketitle 

\begin{abstract} The controllability cost  for the heat equation as the control time $T$ goes to 0 is well-known of the order $e^{C/T}$ for some positive constant $C$,  depending on the controlled domain and for all initial datum. In this paper, we prove that the constant $C$ can be chosen to be  arbitrarily  small if the support of the initial data is sufficiently  close to the controlled domain, but not necessary inside the controlled domain.
The proof is in the spirit  on Lebeau and Robbiano's approach in which a new spectral inequality is established. The main ingredient of the proof of the new spectral inequality is three-sphere inequalities with partial data.   
\end{abstract}



\maketitle

\section{Introduction}

We are interested in the dependence of the  cost of fast controls  for the heat equation on the support (location) of the initial data. 
Let $\omega \subsetneq \Omega$ be a bounded, open subset of $\mR^d$ ($d \ge 1$),  $T> 0$,  $u_0 \in L^2(\Omega)$, and $f \in L^2((0, T) \times \omega )$. Let $A$ be a  Lipschitz, symmetric, uniformly elliptic, matrix-valued function defined in $\Omega$. Consider  the unique solution  $u \in L^2((0, T); H_0^1(\Omega)) \cap C([0, T];  L^2(\Omega))$ of the system 
\begin{equation}\label{CS}
\left\{\begin{array}{cl}
\partial_t u - \dive \big(A(x) \nabla u \big) = f \mathds{1}_{\omega} & \mbox{ in } (0, T) \times \Omega, \\[6pt]
u= 0 & \mbox{ on } (0, T) \times \partial \Omega, \\[6pt]
u(t =0, \cdot) = u_0 & \mbox{ in } \Omega. 
\end{array}\right. 
\end{equation}
Here and in what follows,   $\mathds{1}_{D}$ denotes the characteristic function of a set $D$ of $\mR^d$.  It is well-known from the work of Gilles Lebeau  and Luc Robbiano \cite{LR95},  via spectral inequalities  and the work of Andrei Fursikov and Yu Imanuvilov \cite{FI96}, via Carleman's estimates that one can act on $\omega$ using $f$ to bring $u$ from the initial state $u_0$ (arbitrary) at time $0$ to the final state $0$ at time $T$ (arbitrarily positive).

For $D \subset \Omega$, set 
\begin{equation}\label{def-cost}
c(T, \omega, D) = \mathop{\mathop{\sup}_{\| u_0\|_{L^2(\Omega)} = 1}}_{\supp u_0 \subset D}  \quad \mathop{ \mathop{\inf}_{ \mathop{f \in L^2((0, T) \times   \omega)}_{\mbox{}}}}_{u(T, \cdot) = 0 \; \mathrm{where} \; u \; \mathrm{satisfies} \;  \eqref{CS}}   \| f \|_{L^2((0, T) \times \omega)} .  
\end{equation}
For  $T \in (0, 1)$, one can prove that  
\begin{equation}\label{cost-Omega}
c_1 e^{c_2/ T}  \le c(T, \omega, \Omega) \le C_1 e^{C_2/ T}, 
\end{equation}
for some positive constants $c_1, \, c_2, \, C_1$,  and $C_2$ independent of $T$.  The second inequality follows from the observability inequality \cite{LR95,FI96},  and 
the first inequality was obtained by Luc Miller \cite{M06} and others  \cite{TT07,EZ11}. There is significant literature covering other aspects of the cost of the control for heat equations \cite{CZ00, DE19}, the transport equation with small viscosities \cite{CG05, G10, L12,Lissy15}, and the wave equation \cite{R95, P10, LL19}. The cost of fast controls were also considered for 
linear thermoelasticity \cite{LZ98}, Schr\"odinger equations \cite{Phung01,Miller04,CZ00,Lissy15}, and   plate vibrations \cite{Miller04}. Similar questions were previously addressed in finite dimensions by Thomas Seidman \cite{SeidmanI}.

\medskip 
The goal of this paper is to show a dependence of $c(T, \omega, D)$ on $D$. More precisely, we prove   

\begin{theorem}\label{thm} Let $T \in (0, 1)$  and $\eps > 0$. Assume that  $\omega \Subset \Omega$ is of class $C^2$, and  set, for $r > 0$,  
\begin{equation}\label{def-omega-r}
\omega_r = \Big\{x \in \mR^d; \dist(x, \omega) < r \Big\}. 
\end{equation}
There exist two constants $\delta \in (0, 1)$ and $C_\eps>0$, depending only on $\eps$, $\omega$, $\Omega$,  and the elliptic and Lipschitz constants of $A$, such that  
\begin{equation}
c(T, \omega, \omega_\delta) \le C_\eps e^{\eps/T}. 
\end{equation}
\end{theorem}

\begin{remark} \rm The constants $\delta$ and $C_\eps$ in \Cref{thm} are independent of $T$. 
\end{remark}

When $\omega = \Omega$, the dependence of $c_2$ and $C_2$ on $\Omega$ has been studied extensively, see e.g. \cite{M06,EZ11} and the references therein. Nevertheless,   to our knowledge,  the dependence of the cost on the support of initial datum for the heat equation has not been considered in the literature. \Cref{thm} is new even in one dimensional case. 

\Cref{thm} is expected in the sense that if the support of the initial data is not too far from the controlled region, then it is easier to control. Even in this regard, this intuition is not completely transparent since the propagation speed is infinite and hence the support of the solution at any positive time is generally the whole domain $\Omega$.  Known examples used in the moment method for the heat equations (mainly for one dimensional space) and other equations  give the same size of the control cost for initial datum formed by eigenfunctions of the corresponding operator. From this aspect, \Cref{thm} is thus  unexpected.  

The proof of \Cref{thm} is in the spirit of   Gilles Lebeau and Luc Robbiano's approach \cite{LR95}  in which 
we establish a new spectral inequality.  Let  $0< \lambda_1 \le \lambda_2 \le \dots$ be the sequence of the eigenvalues of the operator $-\dive (A \nabla \cdot)$ with the zero Dirichlet boundary condition,  and let  $e_1, e_2, \dots$ be the corresponding eigenfunctions, i.e., 
\begin{equation}
\left\{\begin{array}{cl}
-\dive(A \nabla e_i) = \lambda_i e_i & \mbox{ in } \Omega, \\[6pt]
e_i = 0 &  \mbox{ on } \partial \Omega. 
\end{array}\right. 
\end{equation}
Assume that $\{ e_i, i \ge 1 \}$ forms an orthogonal basis in $L^2(\Omega)$.  Set, for $\lambda >0$, 
\begin{equation}\label{def-Elambda}
E_{\le \lambda} =  \left\{ \sum_{\lambda_i \le \lambda} a_i e_i(x); a_i \in \mR\right\}. 
\end{equation} 
One of the key elements of Gilles Lebeau and Luc Robbiano's approach is the following spectral inequality 
\begin{equation}\label{SI}
\|v \|_{H^1(\Omega)} \le C e^{C \sqrt{\lambda}} \| v\|_{L^2(\omega)} \quad \forall \,  v \in E_{\lambda}, 
\end{equation} 
where $C$ is a positive constant independent of $\lambda$. 

In this paper, we also follow this approach. Nevertheless, to capture the dependence on the support of the initial datum,  we  use and establish  the following new spectral inequality (compare with \eqref{SI}).

\begin{proposition} \label{pro1} Let  $\eps \in (0, 1)$. There exist two  constants  $\delta \in (0, 1)$ and $C_\eps>0$,  depending only on  $\eps$, $\omega$, $\Omega$, and the elliptic and Lipschitz constants of $A$,  such that, for $\lambda > 0$,  
\begin{equation*}
\| v \|_{L^2(\omega_\delta)} \le C_\eps e^{\eps \sqrt{\lambda}} \| v \|_{L^2(\omega)} \quad \forall v \in E_{\lambda}. 
\end{equation*}
\end{proposition}

\begin{remark}  \rm It is important to emphasize here that the constants $\delta$ and $C_\eps$ in \Cref{pro1} are independent of $\lambda$.  
\end{remark}

\medskip 
The proof of \Cref{pro1} is in the spirit of \cite{LR95}. Nevertheless,  we use three-sphere inequalities with partial data, which was recently established by the author,  to quantitatively capture  the dependence of the support. These inequalities have been derived and applied to the study of cloaking using negative-index materials \cite{Ng-CALR-O,Ng-CALR-O-M}. A typical example of these inequalities is, see \cite[Theorem 2.1]{Ng-CALR-O-M},   

\begin{theorem} \label{thm-3SP}
Let $d \ge 2$, $\Lambda \ge 1$,  $0  < R_1 < R_3 $,   and let $\Gamma = \Big\{x = (x', x_d) \in \partial B_{R_1}; x_d = 0 \Big\}$.  Denote  $O_r = \Big\{x \in \mR^d; \dist(x, \Gamma) < r \Big\}$, $D_r = B_{R_3} \setminus (\overline{B_{R_1} \cup O_r})$, and  $\Sigma_{r} = \partial B_{R_1} \setminus \bar O_{r}$ for $r>0$. For every $\alpha \in (0,  1)$, there exists $r_2 \in (0, R_3 - R_1)$,  depending only on $\alpha$, $\Lambda$, $\Gamma$, 
$R_1$, and $R_3$,   such that for every $r_1 \in (0, r_2)$, there exists $r_0 \in (0, r_1)$,  depending only on $r_1$,  $\alpha$, $\Lambda$, 
$R_1$, and $R_3$,  such that  for  $(d \times d)$ Lipschitz, uniformly elliptic, symmetric, matrix-valued function $\M$ defined in $D_{r_0}$ verifying, in $D_{r_0}$,  
\begin{equation}\label{thm-3SP-cdM}
\Lambda^{-1} |\xi|^2 \le \langle \M (x) \xi, \xi \rangle \le \Lambda |\xi|^2  \; \;  \forall \,  \xi \in \mR^d  \quad \mbox{ and } \quad   |\nabla \M(x) | \le \Lambda, 
\end{equation}
and for  $V  \in [H^1(D_{r_0})]^m$ satisfying  
\begin{equation}\label{fund-thm-Ineq} 
|\dive (\M \nabla V)| \le  \Lambda_1 \big( |\nabla V| + |V|\big)  \mbox{ in } D_{r_0} \mbox{ for some } \Lambda_1 \ge 0, 
\end{equation}
we have 
\begin{equation}\label{thm-3SP-cl}  
\| V \|_{H^1(B_{R_1 + r_2} \setminus B_{R_1 + r_1})} \le C \Big( \| V\|_{H^{1/2}(\Sigma_{r_0})} + \| \M \nabla V \cdot x/|x|\|_{H^{-1/2}(\Sigma_{r_0})} \Big)^\alpha  \|V\|_{H^1(D_{r_0})}^{1 - \alpha}, 
\end{equation}
for some positive constant $C$, depending only on $\alpha$, $\Lambda$,  $\Lambda_1$,  $R_1$, $R_3$,  and $d$. 
\end{theorem}

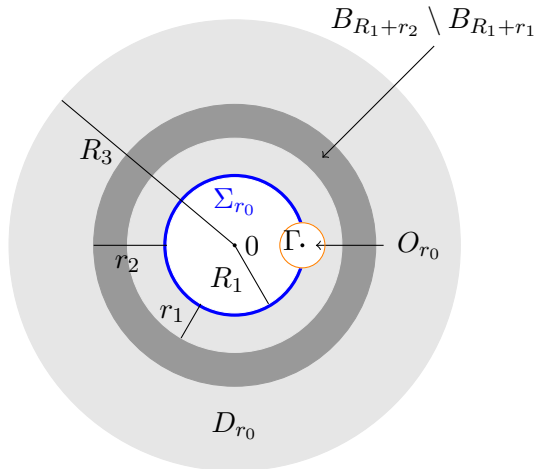
\begin{figure}
\centering
\begin{tikzpicture}[scale=1.5]

\filldraw[gray!20!] (0,0) circle (2);

\filldraw[gray!80!] (0,0) circle (1.25);

\filldraw[gray!20!] (0,0) circle (0.95);

\draw[blue, line width=1mm] (0,0) circle (0.6);

\filldraw[white] (0,0) circle (0.6);

\filldraw[white] (0.6,0) circle (0.2);

\draw[blue] (0,0.6) node[below]{$\Sigma_{r_0}$};

\draw (0,0) node[right]{$0$};

\draw [->] (1.32,0) -- (0.72,0);  

\draw[] (1.35,0) node[right]{$O_{r_0}$};

\draw[violet, orange] (0.6,0) circle (0.2);

\draw (0.68,0.05) node[left]{$\Gamma$};

\fill (0.6,0) circle(0.5pt); 

\fill[black] (0,0) circle(0.5pt); 

\draw (0,0) -- ({0.6*cos(-60)}, {0.6*sin(-60)});

\draw ({0.35*cos(-60)}, {0.35*sin(-60)}) node[left]{$R_1$};

\draw ({0.6*cos(-120)}, {0.6*sin(-120)}) -- ({0.95*cos(-120)}, {0.95*sin(-120)});

\draw ({0.7*cos(-120)}, {0.7*sin(-120)}) node[left]{$r_1$};

\draw ({0.6*cos(-180)}, {0.6*sin(-180)}) -- ({1.25*cos(-180)}, {1.25*sin(-180)});

\draw ({0.95*cos(-180)}, {0.95*sin(-180)}) node[below]{$r_2$};

\draw (0,0) -- ({2*cos(-220)}, {2*sin(-220)});

\draw ({1.6*cos(-220)}, {1.6*sin(-220)}) node[below]{$R_3$};

\draw ({1.6*cos(-90)}, {1.6*sin(-90)}) node[]{$D_{r_0}$};

\draw[->] ({2.5*cos(45)}, {2.5*sin(45)}) -- ({1.1*cos(45)}, {1.1*sin(45)});

\draw[] ({2.5*cos(45)}, {2.5*sin(45)}) node[above]{$B_{R_1 + r_2} \setminus B_{R_1 + r_1}$};

\end{tikzpicture}
\caption{Geometry of \Cref{thm-3SP} in two dimensions}
\label{fig-Thm}
\end{figure}

The geometry of \Cref{thm-3SP} is given in \Cref{fig-Thm}.

\medskip 
We will use a variant of  \Cref{thm-3SP},   given in \Cref{fund-pro},  to derive \Cref{thm}. Nevertheless, we present \Cref{thm-3SP} here to highlight  the difference between the three-sphere inequalities used in this paper and the standard three-sphere ones.  In  \eqref{thm-3SP-cl}, one only uses the information of $\Sigma_{r_0}$ (a portion  of $\partial B_{R_1}$, see \Cref{fig-Thm}) in the first interpolation term. The terminology {\it partial data} comes from this.  The constants $r_1$, $r_2$, and $r_0$ are independent of $\Lambda_1$,  but the constant $C$ does depend on $\Lambda_1$.  If instead of $\Sigma_{r_0}$, one uses $\partial  B_{R_1}$, inequality  \eqref{thm-3SP-cl} 
is then known. Using known three-sphere inequalities and the arguments of the propagation of smallness, one can prove \eqref{thm-3SP} for some $\alpha \in (0, 1)$. Nevertheless, the non-triviality and the novelty of \Cref{thm-3SP} rely on the fact that, for a given arbitrary $\alpha \in (0, 1)$, \eqref{thm-3SP-cl} holds for some $r_0, r_1, r_2$. Even if $v$ is a solution of the Laplace equation in two dimensions, using Hadamard  three-sphere (circles) inequalities  and  the arguments of propagation of smallness, as far as we know, one can only obtain  \eqref{thm-3SP-cl} for some small $\alpha$,  even though one replaces $\Sigma_{r_0}$ by $\partial B_{R_1} \setminus \{x_0 \}$ for some $x_0 \in \partial B_{R_1}$.  The possibility to take $\alpha $ close to $1$  is crucial for the proof of \Cref{thm} where $\eps$ can be arbitrarily small.  This point is also crucial for the cloaking applications considered in \cite{Ng-CALR-O,Ng-CALR-O-M}.  Several applications of \Cref{thm-3SP} concerning variants of Hadamard's three-circle inequalities with partial data are given in \cite{Ng-CALR-O-M}.

\begin{figure}
\centering
\begin{tikzpicture}[scale=1.4]

%\draw[fill=black!20] (-90:1) arc (-90:90:1)-- (90:2) arc (90:-90:2) -- cycle;

\fill[black!10] (-90:2.4) arc (-90:90:2.4)-- (90:0.5) arc (90:-90:0.5) -- cycle;

\fill[black!25] (-90:1.1) arc (-90:90:1.1)-- (90:1.6) arc (90:-90:1.6) -- cycle;

\draw (0,0) -- ({1.1*cos(150-90)},{1.1*sin(150-90)});
\draw ({0.6*cos(150-90)},{0.6*sin(150-90)}) node[right]{$\frac{\hgamma_1 R}{4}$};

\draw (0,0) -- ({1.6*cos(120-110)},{1.6*sin(120-110)});
\draw ({1.4*cos(120-110)},{1.4*sin(120-110)}) node[below]{{$\frac{\hgamma_2 R}{4}$}};

\draw[thick, blue] (0, -2.4) -- (0, -0.5);
\draw[thick, blue] (0, 0.5) -- (0, 2.4);

\draw[->] (-0.5, 2.15) --  (0, 2.15);
\draw[] (-0.5, 2.15) node[left]{$\Sigma_{\gamma_0}$};

\draw [white, domain =-90:90] plot ({0.5*cos(\x)}, {0.5*sin(\x)});

\draw [white, domain =-90:90] plot ({1.6*cos(\x)}, {1.6*sin(\x)});

%\draw [dashed, blue, domain =-90:90] plot ({0.5*cos(\x)}, {0.5*sin(\x)});

\draw[] (0,0) -- ({0.5*cos(30-90)},{0.5*sin(30-90)});
\draw ({0.27*cos(30-90)},{0.27*sin(30-90)}) node[left]{{\small $\hgamma_0 R$}};

\draw ({1.55*cos(75)},{1.55*sin(75)}) node[below]{$Y_{\hgamma_1, \hgamma_2, \frac{R}{4}}$}; 

\draw ({2*cos(-75)},{2*sin(-75)}) node[]{$D_{\gamma_0}$};

\fill[black] (0,0)  circle(0.5pt);

\draw (-0.1,0)  node[left, above]{$0$};

\end{tikzpicture}
\caption{Geometry of $Y_{\gamma_1, \gamma_2, \frac{R}{4}}$, $\Sigma_{\gamma_0}$, and $D_{\gamma_0}$ in two dimensions.}
\label{fig}
\end{figure}
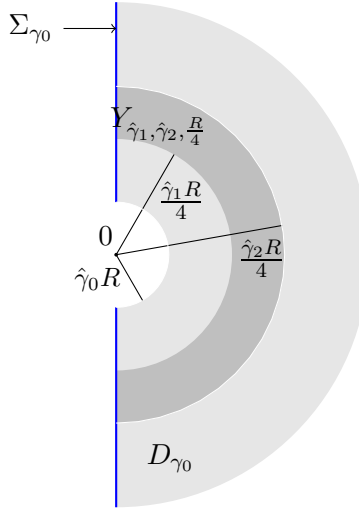

We now introduce some notations to state the local version of \Cref{thm-3SP}, which is used in the proof of \Cref{pro1}.  For $d \ge 2$ and $x = (x_1, x_2, \tx) \in \mR \times \mR \times \mR^{d-2}$, we use the polar coordinate $(\rr, \theta)$ with $\theta \in (-\pi, \pi]$ for the pair $(x_1, x_2)$;  the variable $\tx$ is irrelevant for $d = 2$. 
For $0 < \gamma_1 < \gamma_2 < 1$ and for $R>0$, we denote 
\begin{equation}\label{def-Y-M}
Y_{\gamma_1, \gamma_2,  R} 
= \Big\{ x  \in \mR^d; \; \theta \in (-\pi/2, \pi/2), \;   \gamma_1 R < \rr < \gamma_2 R,  \mbox{ and } |\tx| < R \Big\}, 
\end{equation}
(see also \Cref{fig}). The following variant of \Cref{thm-3SP} 
in a half plane, see \cite[Theorem 3.1]{Ng-CALR-O-M},  is the key ingredient of the proof of \Cref{pro1}.

\begin{proposition} \label{fund-pro}  Let $d \ge 2$,  $\Lambda \ge 1$,  and $R_* < R < R^*$.  Then, for any  $\alpha \in (0, 1)$, there exists a  constant $\hgamma_2 \in (0, 1)$, depending only on $\alpha$, $\Lambda$, $R_*$, $R^*$, and $d$ such that
for every $\hgamma_1 \in (0,  \hgamma_2)$, there exists  $\hgamma_0 \in (0, \hgamma_1)$,  depending only on $\alpha$, $\hgamma_1$,  $\Lambda$, $R_*$, $R^*$, and $d$ such that,   for real, symmetric, uniformly elliptic, Lipschitz matrix-valued functions  $\M$  defined in $D_{\hgamma_0}: = Y_{\hgamma_0, 1,  R}$  verifying, in $D_{\hgamma_0}$, 
\begin{equation}\label{fund-thm-pro-M}
\Lambda^{-1} |\xi|^2 \le \langle \M (x) \xi, \xi \rangle \le \Lambda |\xi|^2  \; \;  \forall \,  \xi \in \mR^d  \quad \mbox{ and } \quad   |\nabla \M(x) | \le \Lambda, 
\end{equation}
and for  $V  \in [H^1(D_{\hgamma_0})]^m$ satisfying
\begin{equation} 
|\dive (\M \nabla V)| \le  \Lambda_1 \big( |\nabla V| + |V|\big)  \mbox{ in } D_{\hgamma_0} \mbox{ for some } \Lambda_1 \ge 0, 
\end{equation}
we have, with $\Sigma_{\hgamma_0} = \partial D_{\hgamma_0} \cap \big\{ x_1= 0 \big\}$,  
\begin{equation}\label{fund-thm-S}
\| V \|_{H^1(Y_{\hgamma_1, \hgamma_2, \frac{R}{4}})} 
 \le C \Big( \|V\|_{H^{1/2}(\Sigma_{\hgamma_0})}+ \|\M \nabla V \cdot \eta_1  \|_{H^{-1/2}(\Sigma_{\hgamma_0})} \Big)^{\alpha}   \| V\|_{H^1(D_{\hgamma_0})}^{1-\alpha},  
\end{equation}
for some positive constant $C$,  depending only on $\alpha, \, \hgamma_1,  \,  \Lambda, \, \Lambda_1, \, R_*, \, R^*$, and  $d$. 
\end{proposition}

Here and in what follows, $\eta_1, \cdots, \eta_d$ denotes the standard basis of $\mR^d$, i.e.,  $\eta_1 = (1,0, \dots, 0)$, \dots, $\eta_d = (0, \dots, 0, 1)$. 

The proof of \Cref{fund-pro} given in \cite{Ng-CALR-O-M}  is quite delicate and involves new (uniform) Carleman's inequalities applied to second-order elliptic equations  in which  the coefficients might be degenerate and in which the geometry of the considered domain  is taken into account in the proof. The proof is much simpler for the case $A = I$ and $d=2$,  but already contains several key ideas \cite{Ng-CALR-O}. 

\medskip 
The paper is organized as follows. \Cref{sect-pro1} is devoted to the proof of \Cref{pro1}. The proof of \Cref{thm} is given in \Cref{sect-thm}.

\section{Spectral inequality}\label{sect-pro1}

This section is devoted to the proof of \Cref{pro1}.  The key ingredient  of the proof is: 

\begin{lemma}\label{lem-3SP} Let  $M$ be a Lipschitz, symmetric, uniformly elliptic,  matrix-valued defined in $\Omega \times (-1, 1)$ and  let $\varphi \in H^1(\Omega \times (-1, 1))$ be such that 
\begin{equation*}
|\dive (M\nabla \varphi)| \le \Lambda (|\nabla \varphi| + |\varphi|) \mbox{ in } \Omega \times (-1, 1). 
\end{equation*}
Set 
$$
\mbox{$D_r = \Big\{X = (x, x_{d+1}) \in \mR^{d+1}; \dist(X , \omega \times \{0\}) < r  \Big\}$ for $r>0$.} 
$$ 
Given $\alpha \in (0, 1)$, there exist two constants $\delta \in (0, 1)$ and $C_\alpha > 0$,  depending only on $\alpha$, $\omega$, $\Omega$, and  the Lipschitz and elliptic constants of $M$,  such that 
\begin{equation}\label{lem-3SP-state}
\|\varphi\|_{H^1(D_\delta)} \le C_\alpha  \Big( \|\varphi \|_{H^{1/2}(\omega \times \{0\})} +  \|M \nabla \varphi \cdot \eta_{d+1} \|_{H^{-1/2}(\omega \times \{0\})} \Big)^\alpha  \| \varphi\|_{H^1(\Omega \times (-1, 1))}^{1 - \alpha}. 
\end{equation}
\end{lemma}

\begin{remark} \rm The constant $\delta$ and $C_\alpha$ are independent of $\varphi$.
\end{remark}

%Here and in what follows, for a subset $\Sigma$ of $\mR^{d+1}$, $H^{-1/2}(\Sigma)$ denotes the dual space of $H^{1/2}_0(\Sigma)$ with the corresponding norm. 

\begin{proof}[Proof of \Cref{lem-3SP}] Since $\omega$ is of class $C^2$, by using local charts and a change of variables, it suffices to prove the following result: Let $\hM$ be a Lipschitz, symmetric, uniformly elliptic,  matrix-valued  function defined in $Q: = (-1, 1)^{d+1}$,  and let $\hvarphi \in H^1(Q)$ be such that 
\begin{equation}\label{lem-3SP-hvarphi}
|\dive (\hM \nabla \hvarphi)| \le \hLambda (|\nabla \hvarphi| + |\hvarphi|) \mbox{ in } Q. 
\end{equation}
Given $\alpha \in (0, 1)$, there exists $\delta > 0$,  depending only on $\alpha$, and  the Lipschitz and elliptic constants of $\hM$,  such that 
\begin{equation}\label{lem-3SP-p0}
\|\hvarphi\|_{H^1(B_\delta)} \le C  
\Big( \|\hvarphi \|_{H^{1/2}(\Sigma)}
 +  \|\hM \nabla \hvarphi \cdot \eta_{d+1} \|_{H^{-1/2} (\Sigma)} \Big)^\alpha \|\hvarphi\|_{H^1(Q)}^{1 - \alpha},
\end{equation}
where $\Sigma  := Q \cap \{x_{d+1} = 0; x_1 < 0\}$  for some positive constant $C$ depending only on $\alpha$, the Lipschitz and elliptic constants of $\hM$, and $\hLambda$. 

 Here and in what follows in this proof, $B_r$ denotes the open ball centered at 0 and of radius $r>0$ in $\mR^{d+1}$. 
 
  It is important to note that in \eqref{lem-3SP-p0}, the norms in the RHS are considered in the set $\Sigma$ which is defined by $Q \cap \{x_{d+1} = 0; x_1 < 0\}$ and is not given by the set 
  $Q \cap \{x_{d+1} = 0\}$.  See a) of \Cref{fig-lem} for the geometry of \eqref{lem-3SP-p0} and b) of \Cref{fig-lem} for the ideas behind using local charts and covering arguments to obtain \eqref{lem-3SP-state} from \eqref{lem-3SP-p0}.

\begin{figure}
\centering
\begin{tikzpicture}[scale=1.8]

\draw [fill=black!5] (-1.2,-1.2) rectangle (1.2,1.2);

\draw [] (-1,1) node{$Q$};

\draw[->] (-1.6,0) -- (1.6,0); 

\draw[] (1.5,0) node[below]{$x_1$}; 

\draw[->] (0,-1.6) -- (0,1.6); 

\draw[] (0,1.5) node[left]{$x_{d+1}$}; 

\draw[red] (0, 0) circle(0.4); 
\draw[<->, dashed] (0, 0) -- ({0.4*cos(30)}, {0.4*sin(30)}); 
\draw[red] (0.1, 0.2) node{$\delta$};

\draw[black!30!green, thick, <->] (-1.2, 0) --(0, 0); 
\draw[black!30!green] (-0.7,0) node[below]{$\Sigma$}; 
 
\draw[] (0, -0.1 - 1.6) node[below]{$a)$}; 
 
\newcommand\x{3.5}
 
\draw [] (-0.8+ \x,1) node{$Q$};
 
\draw[->] (-1.5 + \x,0) -- (1.5 + \x,0); 

\draw[] (1.5 + \x,0) node[below]{$x_1$}; 

\draw[->] (0 + \x,-1.6) -- (0 + \x,1.6); 

\draw[] (0 + \x,1.5) node[left]{$x_{d+1}$}; 

\draw[blue, line width=0.4mm] (-1. + \x, 0) --(1.2 + \x, 0); 

\draw[] (-1 + \x, -1.2) rectangle (1.2 + \x, 1.2); 

\draw[orange, line width=0.8mm] (-0.6 + \x, 0) --(0.8 + \x, 0); 

\draw[dashed] (-0.6 + \x, 0) circle(0.2); 
\draw[dashed] (-0.6 + \x + 0.27, 0) circle(0.2); 
\draw[dashed] (-0.6 + \x + 2*0.27, 0) circle(0.2); 
\draw[dashed] (-0.6 + \x + 3*0.27, 0) circle(0.2); 
\draw[dashed] (-0.6 + \x + 4*0.27, 0) circle(0.2); 
\draw[dashed] (-0.6 + \x + 5*0.27, 0) circle(0.2); 

%\draw[violet] (-0.6 + \x, 0) circle(0.12); 

\draw[violet, thick] (-0.6 + \x, 0.1)  arc(90:270:0.1);

\draw[violet, thick] (0.8 + \x, 0.1)  arc(90:-90:0.1);

\draw[violet, thick] (-0.6 + \x, 0.1)  -- (0.8 + \x, 0.1); 

\draw[violet, thick] (-0.6 + \x, -0.1)  -- (0.8 + \x, -0.1); 

\draw[] (0 + \x, -0.1 - 1.5) node[below]{$b)$};

\end{tikzpicture}
\caption{$a)$: Geometry of inequality  \ref{lem-3SP-p0} in two dimensions with $\Sigma : = Q \cap \{x_{d+1} = 0; x_1 < 0\}$. $b)$ The way to obtain \eqref{lem-3SP-state} from \eqref{lem-3SP-p0} for $d=1$; $\omega$ is the orange interval, $\Omega$ is the blue interval, $D_\delta$ is the region whose boundary is violet.}
\label{fig-lem}
\end{figure}
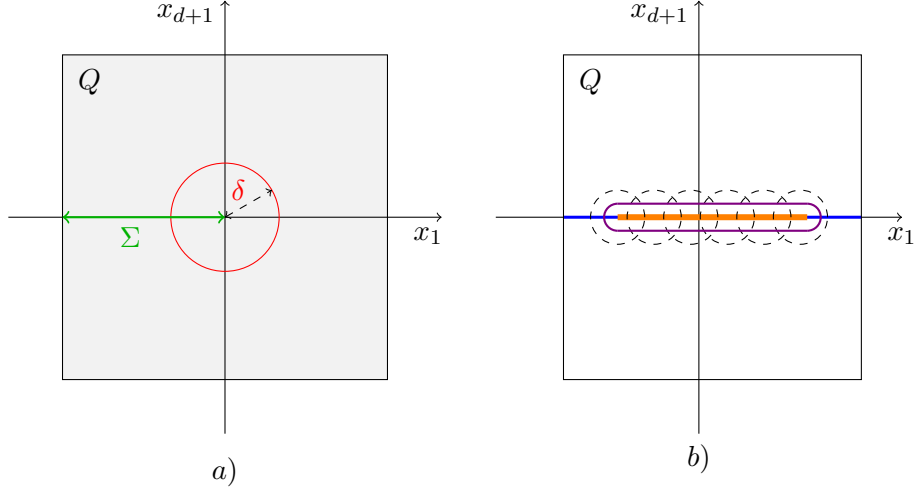

We will make a change of variables  in order to apply \Cref{fund-pro}. To this end, for $X = (x_1, \cdots, x_{d+1}) \in Q \setminus \big\{x_{d+1} = 0; x_1 \le 0\big\}$, define
$$
\cR(X) = (y_1, x_2, \cdots, x_d, y_{d+1}), 
$$
with $(y_1, y_{d+1}) = \hat r  e^{i \theta /2}$ if $(x_1, x_{d+1}) = \hat r e^{i \theta}$ for $\hat r  >0$,  and $\theta \in (-\pi, \pi)$.

Set 
\begin{equation*}
\hvarphi_1 = \hvarphi \circ \cR^{-1} \mbox{ in } \hat Q_1: =  \cR \Big(Q \setminus \big\{x_{d+1} = 0; x_1 \le  0 \big\} \Big). 
\end{equation*}
Set 
$$
f (x) = \dive (\hM \nabla \hvarphi) (x) \mbox{ in } Q,  \quad  f_1 (x) =  \frac{f}{|\det (\nabla \cR)|} \circ \cR^{-1} (x) \mbox{ in } Q_1, 
$$
and 
\begin{equation}\label{lem-3SP-M1}
 \hM_1 = \frac{\nabla \cR \hM \nabla \cR \tr}{|\det (\nabla \cR)|} \circ \cR^{-1} \mbox{ in } \hat Q_1.
\end{equation}
It is clear from \eqref{lem-3SP-M1} that the elliptic and Lipschitz constants of $\hM_1$ are bounded by the elliptic and Lipschitz constants of $\hM$,  up to a constant $C$,  depending only on $d$.   

Since $\dive (\hM \nabla \hvarphi)  = f$ in $Q$, it follows from  a change of variables that 
\begin{equation}\label{lem-3SP-p1}
\dive (\hM_1 \nabla \hvarphi_1) = f_1 \mbox{ in } Q_1. 
\end{equation}
We have
$$
 \det (\nabla \cR) (x) = 1/2 \mbox{ for } x \in Q, 
$$
$$
|\nabla \hvarphi (x)| \le  |\nabla \cR (x)| |\nabla \hvarphi_1 \circ \cR (x)| \le C |\nabla \hvarphi_1 \circ \cR (x)| \mbox{ for } x \in Q.   
$$
Since $|f| \le C \big( |\nabla \varphi| + |\varphi| \big)$  in $Q$ by \eqref{lem-3SP-hvarphi}, we derive from  \eqref{lem-3SP-p1} that 
\begin{equation*}
|\dive (\hM_1 \nabla \hvarphi_1)| \le \hLambda_1  \big( |\nabla \varphi_1| + |\varphi_1| \big) \mbox{ in } Q_1, 
\end{equation*}
for some $\hLambda_1 > 0$,  depending only on $\Lambda$ and $d$.  

Set
\begin{equation}
\Gamma_{1, +} = \Big\{(x_1, \dots, x_{d+1}); x_{1} = 0, \;  x_{d+1}  >  0, \;  x_j \in (-1, 1) \mbox{ for } 2 \le j \le d  \Big\}
\end{equation}
and
\begin{equation}
\Gamma_{1, -} =  \Big\{(x_1, \dots, x_{d+1}); x_{1} = 0, \;  x_{d+1} < 0, \;  x_j \in (-1, 1) \mbox{ for } 2 \le j \le d  \Big\}.
\end{equation}
Apply \Cref{fund-pro} to $\hvarphi_1$ with $R = 1$, $\hgamma_1 = \hgamma_2/2$,  and $\hgamma_0 = 0$, and in $\mR^{d+1}$ with $(x_1, x_2, \tx)$ being replaced by $\big(x_1, x_{d+1}, (x_2, \cdots, x_d) \big)$. There exists $\hgamma_2 > 0$ such that 
\begin{multline*}
\| \hvarphi_1\|_{H^1\big( (B_{\hgamma_2} \setminus B_{\hgamma_2/2}) \cap \{x_1 > 0\} \big)} \le C 
\| \hvarphi_1\|_{H^1( \hat Q_1)}^{1 - \alpha}\\[6pt]
\times \Big( \| \hvarphi_1 \|_{H^{1/2} \big(\Gamma_{1, +}  \big)} + \| \hM_1 \nabla \hvarphi_1 \cdot \eta_{1} \|_{H^{1/2} \big( \Gamma_{1, +} \big)} + \| \hvarphi_1 \|_{H^{1/2} \big(\Gamma_{1, -}  \big)} + \| \hM_1 \nabla \hvarphi_1 \cdot \eta_{1} \|_{H^{1/2} \big( \Gamma_{1, -} \big)} \Big)^\alpha. 
\end{multline*}
Since $\hvarphi = \hvarphi_1 \circ \cR$ in $Q$, it follows from a change of variables, see e.g. \cite[Lemma 2]{Ng-Complementary}, that 
\begin{equation*}
\|\hvarphi\|_{H^1\big((B_{\hgamma_2} \setminus B_{\hgamma_2/2}) \setminus \{x_{d+1} = 0; x_1 < 0 \} \big)} \le C  \|\hvarphi\|_{H^1(Q)}^{1 - \alpha}   \Big( \|\hvarphi \|_{H^{1/2}(\Sigma)} +  \|\hM \nabla \hvarphi \cdot \eta_{d+1} \|_{H^{-1/2}(\Sigma)} \Big)^\alpha.  
\end{equation*}
Since $\hvarphi \in H^1(Q)$,  and hence in particular $\hvarphi \in H^1(B_{\hgamma_2})$, we obtain  
\begin{equation}\label{lem-3SP-cl1}
\|\hvarphi\|_{H^1(B_{\hgamma_2} \setminus B_{\hgamma_2/2})} \le C  \|\hvarphi\|_{H^1(Q)}^{1 - \alpha} \Big( \|\hvarphi \|_{H^{1/2}(\Sigma)} +  \|\hM \nabla \hvarphi \cdot \eta_{d+1} \|_{H^{-1/2}(\Sigma)} \Big)^\alpha. 
\end{equation}
Using the fact 
\begin{equation*}
|\dive (\hM\nabla \hvarphi)| \le \hLambda (|\nabla \varphi| + |\varphi|) \mbox{ in } Q, 
\end{equation*}
and $\hM$ is symmetric, uniformly elliptic and  Lipschitz, one has \footnote{One can prove \eqref{lem-3SP-cl2} using a contradiction argument and the unique continuation principle.} 
\begin{equation}\label{lem-3SP-cl2}
\|\hvarphi \|_{H^1(B_{\gamma_2})}  \le C \|\hvarphi\|_{H^1(B_{\gamma_2} \setminus B_{\gamma_2/2})}.  
\end{equation}
Assertion \eqref{lem-3SP-p0} now follows from \eqref{lem-3SP-cl1} and \eqref{lem-3SP-cl2} with $\delta = \gamma_2$. The proof is complete. 
\end{proof}

\begin{remark}\rm One of the key points of the proof is the assertion \eqref{lem-3SP-p0}. This assertion is known if one replaces the set $Q \cap \{x_{d+1} = 0; x_1 < 0\}$ by $Q \cap \{x_{d+1} = 0\}$ and the proof in this case can be done as in \cite{LR95}. However, this does not imply \eqref{lem-3SP-state}.  The proof of \eqref{lem-3SP-p0} follows from  \Cref{fund-pro}, which  is non-trivial. 
\end{remark}

We are ready to give 

\begin{proof}[Proof of \Cref{pro1}] Since $v \in E_{\le \lambda}$, there exists $a_i \in \mR$ with $\lambda_i \le \lambda$,  such that 
$$
v (x) = \sum_{\lambda_i \le  \lambda} a_i e_i  (x) \mbox{ in } \Omega. 
$$
As in the spirit of \cite{LR95}, set, with $X = (x, x_{d+1}) \in \Omega \times \mR$, 
\begin{equation*}
V(X) = \sum_{\lambda_i \le  \lambda}  \lambda_i^{-1/2}a_i \sinh(\lambda_i^{1/2} x_{d+1}) e_i (x),    
\end{equation*}
where $\sinh t= \frac{1}{2} (e^{t} - e^{-t})$ for $t \in \mR$. Since $-\dive_x (A(x) \nabla e_i(x)) = \lambda_i e_i (x)$ in $\Omega$, it follows that 
\begin{equation}\label{pro1-sys-V}
\left\{\begin{array}{c}
\partial_{x_{d+1}}^2 V + \dive_x \big(A(x) \nabla_x V\big)  = 0 \mbox{ in } \Omega \times \mR, \\[6pt]
V(X) = 0 \mbox{ for } X \in \Omega \times \{0 \}, \\[6pt]
\partial_{x_{d+1}} V(X) = v(x) \mbox{ for } X \in \Omega \times \{0 \}. 
\end{array}\right.
\end{equation}

Given $\alpha \in (0,1)$, by applying  \Cref{lem-3SP} to $V$, there exist two  constants $\delta  = \delta(\alpha) \in (0, 1)$ and $C_\alpha > 0$,  depending only on $\alpha$,  $\omega$,  $\Omega$, and the elliptic and Lipschitz constants of $A$,  such that 
\begin{equation}\label{pro1-p0-1}
\| V \|_{H^1(D_{2\delta})} \le C_\alpha \Big( \| V \|_{H^{1/2}(\omega \times \{0 \})} + \| \partial_{x_{d+1}} V \|_{H^{-1/2}(\omega \times \{0 \})} \Big)^{\alpha}  \|V \|_{H^1(\Omega \times (-1, 1))}^{1 - \alpha}. 
\end{equation}

Using \eqref{pro1-sys-V}, we derive from  \eqref{pro1-p0-1} that 
\begin{equation}\label{pro1-p0}
\| V \|_{H^1(D_{2\delta})} \le C_\alpha \| v \|_{L^2(\omega)}^{\alpha}  \|V \|_{H^1(\Omega \times (-1, 1))}^{1 - \alpha}. 
\end{equation}
Since $A$ is Lipschitz, by the regularity theory of elliptic equations \footnote{One can directly apply the quotient method due to Louis Nirenberg \cite{Nirenberg59}.}, one has 
\begin{equation*}
 \| \partial_{x_{d+1}} V \|_{L^2(\omega_{\delta})} \le C_\alpha \| V \|_{H^1(D_{2\delta})}, 
\end{equation*}
which yields 
\begin{equation}\label{pro1-p1}
\| v \|_{L^2(\omega_{\delta})} \le C_\alpha \| V \|_{H^1(D_{2\delta})}. 
\end{equation}

On the other hand, by the standard spectral inequality \eqref{SI}, one gets
 \begin{equation}\label{pro1-p2}
 \|V \|_{H^1(\Omega \times (-1, 1))} \le C e^{C \sqrt{\lambda}} \| v\|_{L^2(\omega)}, 
\end{equation}
for some positive constant $C$,  depending only on $\omega$, $\Omega$, and the elliptic and Lipschitz constants of $A$.

Combining \eqref{pro1-p0}, \eqref{pro1-p1} and \eqref{pro1-p2} yields 
\begin{equation}\label{pro1-p3}
\| v \|_{L^2(\omega_\delta)}
 \mathop{\le}^{\eqref{pro1-p1}}  
 C_\alpha \| V \|_{H^1(D_{2\delta})}
  \mathop{\le}^{\eqref{pro1-p0}} 
C_\alpha \| v \|_{L^2(\omega)}^{\alpha}  \|V \|_{H^1(\Omega \times (-1, 1))}^{1 - \alpha} 
\mathop{ \le }^{\eqref{pro1-p2}} C_\alpha e^{C (1 - \alpha) \sqrt{\lambda}} \| v\|_{L^2(\omega)}. 
\end{equation}
By choosing $\alpha$ such that $C (1 - \alpha) = \eps$,  we derive from \eqref{pro1-p3} that 
\begin{equation}
\| v \|_{L^2(\omega_\delta)} \le C_\eps e^{\eps \sqrt{\lambda}} \| v\|_{L^2(\omega_\delta)}. 
\end{equation}
The proof is complete. 
\end{proof}

\section{Proof of \Cref{thm}} \label{sect-thm}

The proof of \Cref{thm} is based on the following lemma,  which will be derived from the spectral inequality stated in \Cref{pro1}. 

\begin{lemma}\label{lem2} Let $0 < T < 1$, $\lambda > 0$,  and let $v_0 \in E_\lambda$. Let $v \in L^2((0, T), H_0^1(\Omega)) \cap C([0, T], L^2(\Omega))$ be the unique solution of the system 
\begin{equation}\label{lem2-eqv}
\left\{\begin{array}{cl}
\partial_t v - \dive(A \nabla v) = 0 & \mbox{ in } (0, T) \times \Omega, \\[6pt]
v= 0 & \mbox{ on } (0, T) \times \partial \Omega, \\[6pt]
v(t =0, \cdot) = v_0 & \mbox{ in } \Omega. 
\end{array}\right. 
\end{equation}
For $\eps > 0$, there exist two constants  $\delta \in (0, 1)$ and $C_\eps > 0$,  depending only  on $\eps$, $\omega$, and $\Omega$, and the elliptic and Lipschitz constants of $A$,   such that 
\begin{equation*}
\| v(T, \cdot)\|_{L^2(\omega_\delta)} \le C_\eps \delta^{-1} T^{-1/2}e^{\eps \sqrt{\lambda}} \| v\|_{L^2((0, T) \times \omega)}. 
\end{equation*}
\end{lemma}

Recall that $\omega_r$ is defined in \eqref{def-omega-r}. 

\begin{remark} \rm The constants $\delta$ and $C_\eps$ in \Cref{lem2} are independent of $\lambda$ and $T$. 
\end{remark}

\begin{proof} By \Cref{pro1}, there exist $\delta \in (0, 1)$ and $C_\eps > 0$,  such that 
\begin{equation}\label{lem2-p1-*}
\| \xi \|_{L^2(\omega_{2 \delta})} \le C_\eps  e^{\eps \sqrt{\lambda}} \| \xi\|_{L^2(\omega)} \mbox{ for } \xi \in E_{\le \lambda}. 
\end{equation}
Since  $v_0 \in E_{\le \lambda}$, it follows that $v(t, \cdot) \in E_{\le \lambda}$ for $t \in (0, T)$. We derive from \eqref{lem2-p1-*} that 
\begin{equation}\label{lem2-p1}
\| v(t, \cdot) \|_{L^2(\omega_{2 \delta})} \le C_\eps  e^{\eps \sqrt{\lambda}} \| v(t, \cdot)\|_{L^2(\omega)}  \mbox{ for } t \in (0, T). 
\end{equation}

Fix $\varphi \in C^\infty_c(\mR^d)$ such that $0 \le \varphi \le 1$,  $\varphi = 1$ in $\omega_\delta$, $\supp \varphi \subset \omega_{2\delta}$, and  $|\nabla_x^\alpha \varphi| \le C / \delta^{|\alpha|}$ for all multi-indices  $\alpha$ with $|\alpha| \le 2$.  Here and in what follows in this proof, $C$ denotes a positive constant, depending only on $\omega$, $\Omega$, and the elliptic and Lipschitz constants of $A$. 

Set 
\begin{equation}\label{lem2-def-u}
u (t, x) = \varphi (x) v (t, x) \mbox{ in } (0, T) \times \Omega, 
\end{equation}
and denote 
\begin{equation}\label{lem2-def-g}
g(t, x) = - \Big( 2 \langle A(x)  \nabla v (t, x), \nabla \varphi(x)  \rangle +  v (t, x) \dive(A (x) \nabla \varphi (x)) \Big) \mbox{ in } (0, T) \times \Omega, 
\end{equation}
where $\langle \cdot, \cdot \rangle$ denotes the standard scalar product in $\mR^d$. 

We derive from \eqref{lem2-eqv} and the symmetry of $A$ that  
\begin{equation*}
\left\{\begin{array}{cl}
\partial_t u - \dive (A \nabla  u) =  g  & \mbox{ in } (0, T) \times \Omega, \\[6pt]
u= 0 & \mbox{ on } (0, T) \times \partial \Omega.  
\end{array}\right. 
\end{equation*}
  Multiplying the equation of $u$ by $u$ and integrating by parts in $(t, T) \times \Omega$, we obtain, for $0 \le t \le T$, 
\begin{multline}\label{lem2-p2}
\frac{1}{2}\int_{\Omega} |u(T, x)|^2 \, dx  + \int_t^T \int_{\Omega} \langle A(x)\nabla u(s, x), \nabla u (s, x) \rangle \, dx \, ds \\[6pt] 
= \frac{1}{2} \int_{\Omega} |u(t, x)|^2 \, dx  +   \int_t^T \int_{\Omega} g (s, x) u (s, x) \, dx  \, ds.   
\end{multline}

We next estimate the last term of \eqref{lem2-p2}. Since, for $x \in \Omega$ and $s \in (0, T)$, 
\begin{multline*}
\langle A (x) \nabla v (s, x),  \nabla \varphi (x) \rangle  u (s, x)   \mathop{=}^{\eqref{lem2-def-u}}  \langle A (x) \varphi (x) \nabla v (s, x),  \nabla \varphi (x) \rangle  v (s, x) \\[6pt]
 \mathop{=}^{\eqref{lem2-def-u}}  \langle A(x) (\nabla u(s, x) -  v(s, x) \nabla \varphi(x) ), \nabla \varphi (x) \rangle  v(s, x), 
\end{multline*}
it follows that, for $s \in (0, T)$,  
\begin{equation}\label{lem2-p3}
\left|\int_{\Omega} \langle A (x) \nabla v (s, x),  \nabla \varphi (x) \rangle  u (s, x) \, dx \right| \le C \int_{\Omega} \delta^{-1} |\nabla u(s, x)| |v(s, x)| \, dx  + \int_{\Omega} \delta^{-2}  |v(s, x)|^2 \, dx. 
\end{equation}
We also have,    for $x \in \Omega$ and $s \in (0, T)$, 
\begin{equation}\label{lem2-p4}
|v(s, x) \dive(A (x) \nabla \varphi (x)) | |u(s, x)| \le C \delta^{-2} |v(s, x)|^2, 
\end{equation}
since $A$ is Lipschitz and $|\nabla_x^\alpha \varphi| \le C / \delta^{|\alpha|}$ for all multi-indices  $\alpha$ with $|\alpha| \le 2$. Combining \eqref{lem2-p3} and  \eqref{lem2-p4}  yields
\begin{multline}\label{lem2-p5}
 \left|\int_t^T \int_{\Omega} g (s, x) u (s, x) \, dx  \, ds \right| \le C \int_t^T \int_{\Omega} \delta^{-1} |\nabla u(s, x)| |v(s, x)|  \, dx \, ds  \\[6pt]
 + \int_t^T \int_{\Omega} \delta^{-2}  |v(s, x)|^2 \, dx \, ds. 
\end{multline}

Using \eqref{lem2-p5} and the ellipticity of $A$,  and applying Young's inequality, we derive from \eqref{lem2-p2}  that, for $t \in (0, T)$,  
\begin{equation*}
\int_{\Omega} |u(T, x)|^2 \, dx \le \int_{\Omega} |u(t, x)|^2  \, dx +  C \delta^{-2} \int_t^T \int_{\Omega} |v(s, x)|^2 \, dx \, ds. 
\end{equation*}
Integrating the above inequality with respect to $t$ from $0$ to $T$,  we derive that 
\begin{equation*}
\int_{\Omega} |u(T, x)|^2 \, dx  \le C \delta^{-2} T^{-1} \int_{0}^T \int_{\Omega} |v(s, x)|^2 \, dx \, ds.  
\end{equation*}
Since $v = \varphi u$,  $0 \le \varphi \le 1$,  $\varphi = 1$ in $\omega_\delta$, and $\supp \varphi \subset \omega_{2 \delta}$, it follows that 
\begin{equation*}
\int_{\omega_\delta} |v(T, x)|^2 \, dx \le C \delta^{-2} T^{-1} \int_{0}^T \int_{\omega_{2\delta}} |v(t, x)|^2 \, dx \, dt. 
\end{equation*}
We derive from \eqref{lem2-p1} that 
\begin{equation*}
\int_{\omega_\delta} |v(T, x)|^2 \, dx \le  C_\eps \delta^{-2} T^{-1} e^{\eps \sqrt{\lambda}} \int_{0}^T \int_{\omega} |v(t, x)|^2 \, dx \, dt, 
\end{equation*}
which is the conclusion. The proof is complete. 
\end{proof}

We are ready to give 

\begin{proof}[Proof of \Cref{thm}]
Fix $\lambda = c_0 / T^2$ where $c_0$ is a large positive constant determined later. Set 
\begin{equation}\label{def-H}
H := \left\{ \sum_{\lambda_i \le \lambda} a_i e^{-\lambda_i (T/3 - t)} e_i (x) ; a_i \in \mR, x \in \Omega, t \in (0, T/3) \right\} \subset L^2\big((0, T/3) \times \Omega\big). 
\end{equation}
Equip $H$ with the standard scalar product in $L^2\big((0, T/3) \times \Omega \big)$. Then,  $H$ is a Hilbert space (of finite dimensions). 

Let $\varphi \in H$, and set 
\begin{equation}\label{def-v}
v(t, x) = \varphi(T/3 - t, x) \mbox{ for } (t, x) \in (0, T/3) \times \Omega. 
\end{equation}
It follows from the definition of $H$ in \eqref{def-H} that 
\begin{equation*}
\left\{\begin{array}{cl}
\partial_t v - \dive(A \nabla v) = 0 & \mbox{ in } (0, T/3) \times \Omega, \\[6pt]
v= 0 & \mbox{ on } (0, T/3) \times \partial \Omega, 
\end{array}\right. 
\end{equation*}
and moreover, $v(t = 0, \cdot) \in E_{\le \lambda}$.

By \Cref{lem2},  there exist two constants  $\delta \in (0, 1)$ and $C_\eps > 0$,  depending only on $\eps$, $\omega$, $\Omega$, $c_0$, and 
the Lipschitz and elliptic constants of $A$,  such that
$$
\| v(T/3, \cdot) \|_{L^2(\omega_\delta)} \le C_\eps \delta^{-1} T^{-1/2} \| v \|_{L^2 \big( (0, T/3) \times \Omega \big)}. 
$$
This implies, by \eqref{def-v},   
\begin{equation}\label{thm-p1}
\|\varphi(0, \cdot)\|_{L^2(\omega_\delta)} \le C_\eps  \delta^{-1} T^{-1/2} e^{\eps/ T } \|  \varphi\|_{L^2\big((0, T/3) \times \omega \big)}. 
\end{equation}
Fix such constants $\delta$ and $C_\eps$. 

Fix $u_0 \in L^2(\Omega)$ with $\supp u_0 \subset \omega_\delta$. We will construct a control with support in $(0, T) \times \omega$,  which steers $u_0$ from time 0 to $0$ at time $T$ for which  the cost is bounded by $C_\eps e^{\eps/ T} \| u_0\|_{L^2(\Omega)}$. 

Since $u_0 \in L^2(\Omega)$ with $\supp u_0 \subset \omega_\delta$, using the  Riesz  representation theorem,  we derive from \eqref{thm-p1} that  there exists $f_1 \in H$,  such that 
\begin{equation}\label{thm-f1}
\int_{\Omega} u_0 (x) \varphi(0, x) \, dx  = \int_{0}^{T/3} \int_{\omega} f_1 (s, x) \varphi (s, x) \, dx \, ds \mbox{ for } \varphi \in H, 
\end{equation}
and 
\begin{equation}\label{thm-norm-f1} 
\| f_1 \|_{L^2\big((0, T/3) \times \omega \big)} \le C_\eps \delta^{-1} T^{-1/2} e^{\eps/ T} \| u_0\|_{L^2(\Omega)}. 
\end{equation}

Let $u_1 \in L^2((0, T/3); H_0^1(\Omega)) \cap C([0, T/3];  L^2(\Omega))$ be the unique solution of the system 
\begin{equation}\label{thm-sys-u1}
\left\{\begin{array}{cl}
\partial_t u_1 - \dive \big(A(x) \nabla u_1 \big) = f_1\mathds{1}_\omega & \mbox{ in } (0, T/3) \times \Omega, \\[6pt]
u_1= 0 & \mbox{ on } (0, T/3) \times \partial \Omega, \\[6pt]
u_1(t =0, \cdot) = u_0 & \mbox{ in } \Omega. 
\end{array}\right. 
\end{equation}
Since 
$$
\left\{\begin{array}{cl}
\partial_t \varphi + \dive(A \nabla \varphi) = 0 \mbox{ in } (0, T/3) \times \Omega, \\[6pt]
\varphi = 0 \mbox{ on } (0, T/3) \times \partial \Omega, 
\end{array}\right.
 \mbox{ for } \varphi \in H, 
$$
multiplying the equation of $u_1$ by $\varphi$ ($\in H$) and integrating by parts in $(0, T/2) \times \Omega$, we obtain 
\begin{equation*}
\int_{\Omega} u_1 (T/ 3, x) \varphi(T/ 3, x) \, dx - \int_{\Omega} u_1 (0, x) \varphi(0, x) \, dx   =  \int_0^{T/3} \int_{\omega} f_1(s, x) \varphi(s, x) \, dx \, ds  \mbox{ for } \varphi \in H.  
\end{equation*}
Using \eqref{thm-f1}, we derive that 
\begin{equation}\label{thm-proj1}
\int_{\Omega} u_1 (T/ 3, x) \varphi(T/ 3, x) \, dx = 0   \mbox{ for } \varphi \in H.  
\end{equation}
In other words, the projection of $u(T/3, \cdot)$ into $E_{\le \lambda}$ is 0. Thus,   
\begin{equation}\label{thm-proj-u1}
u_1(T/3, x) = \sum_{\lambda_i > \lambda} \langle u_1(T/ 3, \cdot),  e_i \rangle_{L^2(\Omega)}  e_i (x) \mbox{ in } \Omega, 
\end{equation}
where $\langle \cdot, \cdot \rangle_{L^2(\Omega)}$ denotes the standard scalar product in $L^2(\Omega)$.

On the other hand, by the standard energy estimate, we have 
\begin{equation*}
\int_{\Omega} |u_1(T/3, x)|^2 \, dx \le 2 \int_{\Omega} |u_1(0, x)|^2 \, dx + C \int_0^{T/3} \int_{\omega} |f_1(s, x)|^2 \, ds \, dx. 
\end{equation*} 
We derive from \eqref{thm-norm-f1} that 
\begin{equation}\label{thm-u1-T/3}
\| u_1(T/3, \cdot) \|_{L^2(\Omega)} \le C_\eps \delta^{-1} T^{-1/2} e^{\eps / T} \| u_0\|_{L^2(\Omega)}. 
\end{equation} 

Let $u_2 \in L^2((T/3, 2T/3); H_0^1(\Omega)) \cap C([T/3, 2T/3];  L^2(\Omega))$ be the unique solution of the system 
\begin{equation}\label{thm-sys-u2}
\left\{\begin{array}{cl}
\partial_t u_2 - \dive \big(A(x) \nabla u_2 \big) = 0 & \mbox{ in } (0, T/3) \times \Omega, \\[6pt]
u_2= 0 & \mbox{ on } (0, T/3) \times \partial \Omega, \\[6pt]
u_2(t =T/3, \cdot) = u_1(T/3, \cdot) & \mbox{ in } \Omega. 
\end{array}\right. 
\end{equation} 
It follows from \eqref{thm-proj-u1} and \eqref{thm-u1-T/3} that 
\begin{equation*}
\| u_2 (2T/3, \cdot) \|_{L^2(\Omega)} \le e^{-\lambda T/3} \| u_2 (T/3, \cdot) \|_{L^2(\Omega)} \le C_\eps \delta^{-1} T^{-1/2} e^{\eps / T -\lambda T/3} \| u_0\|_{L^2(\Omega)}, 
\end{equation*}
which yields, since $\lambda = c_0 / T^2$, 
\begin{equation}\label{thm-u2-norm}
\| u_2 (2T/3, \cdot) \|_{L^2(\Omega)} \le  C_\eps \delta^{-1} T^{-1/2} e^{\eps / T - c_0/ (3T)} \| u_0\|_{L^2(\Omega)}. 
\end{equation}

On the other hand, there exists $f_3 \in L^2((2T/3, T) \times \Omega)$ with support in $[2T/3, T] \times \omega$, such that 
\begin{equation}\label{thm-norm-f3}
\| f_3 \|_{L^2\big((2T/3, T) \times \Omega \big)} \le C e^{C/T} \| u_2(2T/3, \cdot)\|_{L^2(\Omega)}, 
\end{equation}
and 
\begin{equation}\label{thm-cd-u3}
u_3(T, \cdot) = 0 \mbox{ in } \Omega, 
\end{equation}
where  $u_3 \in L^2\big((2T/3, T); H_0^1(\Omega) \big) \cap C\big([2T/3, T];  L^2(\Omega) \big)$ is the unique solution of the system
\begin{equation}\label{thm-sys-u3}
\left\{\begin{array}{cl}
\partial_t u_3 - \dive \big(A(x) \nabla u_3 \big) = f_3 & \mbox{ in } (2T/3, T) \times \Omega, \\[6pt]
u_3= 0 & \mbox{ on } (T/2, T) \times \partial \Omega, \\[6pt]
u_3(t =2T/3, \cdot) = u_2(2T/3, \cdot) & \mbox{ in } \Omega. 
\end{array}\right. 
\end{equation}

Define $f \in L^2 \big( (0, T) \times \Omega \big)$ as follows 
\begin{equation}
f(t, x)  = \left\{ \begin{array}{cl} f_1 \mathds{1}_\omega & \mbox{ in } (0, T/3) \times \Omega, \\[6pt]
0 & \mbox{ in } (T/3, 2T/3) \times \Omega, \\[6pt]
 f_3 & \mbox{ in } (2T/3, T) \times \Omega. 
\end{array}\right.
\end{equation}
Since $\supp f_3 \subset [2T/3, T] \times \omega$, it follows that 
$$
\supp f \subset [0, T] \times \bar \omega. 
$$
Let $u \in L^2\big((0, T); H_0^1(\Omega) \big) \cap C\big([0, T];  L^2(\Omega) \big)$ be the unique solution of the system
\begin{equation}\label{thm-sys-u}
\left\{\begin{array}{cl}
\partial_t u - \dive \big(A(x) \nabla u \big) = f  & \mbox{ in } (0, T) \times \Omega, \\[6pt]
u= 0 & \mbox{ on } (T/2, T) \times \partial \Omega, \\[6pt]
u(t =0, \cdot) = u_0 & \mbox{ in } \Omega. 
\end{array}\right. 
\end{equation}
It follows from \eqref{thm-sys-u1}, \eqref{thm-sys-u2}, \eqref{thm-cd-u3}, and \eqref{thm-sys-u3} that 
\begin{equation*}
u_3(T, \cdot) = 0 \mbox{ in } \Omega. 
\end{equation*}
Combining \eqref{thm-norm-f1} and \eqref{thm-norm-f3}, and using \eqref{thm-u2-norm}, we deduce that  
\begin{equation*}
\| f \|_{L^2((0, T) \times \Omega)} \le C_\eps \delta^{-1} T^{-1/2} e^{\eps/ T} \| u_0\|_{L^2(\Omega)} \left(1 + e^{-c_0/(3T) + C/ T}  \right).  
\end{equation*}
By fixing $c_0$ such that $c_0 / 3 \ge C$, we obtain 
\begin{equation*}
\| f \|_{L^2((0, T) \times \Omega)} \le C_\eps \delta^{-1} T^{-1/2} e^{\eps/ T} \| u_0\|_{L^2(\Omega)}.  
\end{equation*}
The conclusion  follows by replacing $\eps$ by $\eps/2$ and noting that 
$$
T^{-1/2} e^{\eps/ (2T) } \le C_\eps e^{\eps/ T};  
$$
this follows by considering the case $T \ge \eps$ and the case $0< T < \eps$.  
\end{proof}

\begin{remark} \rm The conclusion of \Cref{thm} also holds if in the definition of $c(T,  w, D)$, one additionally requires that $\supp f \Subset [0, T] \times \omega$.  The conclusion in this case follows by applying the established result for the set $\big \{x \in \omega; \dist (x, \partial \omega) \ge \gamma \big\}$ for small $\gamma$ after noting that the constant $\delta$ for such a set is independent of $\gamma$ for small $\gamma$.  
\end{remark}

\medskip 
\noindent{\bf Acknowledgement:} The author thanks Jean-Michel Coron for his interest in the problem and for many interesting discussions. The author also thanks Kim Dang Phung for discussions on the approach of Gilles Lebeau and Luc Robbiano. This work was completed during his visit to Laboratoire Jacques Louis Lions. The author thanks the laboratory  for its hospitality.

\providecommand{\bysame}{\leavevmode\hbox to3em{\hrulefill}\thinspace}
\providecommand{\MR}{\relax\ifhmode\unskip\space\fi MR }
\providecommand{\MRhref}[2]{%
  \href{http://www.ams.org/mathscinet-getitem?mr=#1}{#2}
}
\providecommand{\href}[2]{#2}

\end{document}